\documentclass[12pt, reqno]{amsart}
\usepackage{color}
\usepackage{amsmath}
\usepackage{amssymb}
\usepackage{amsfonts}
\usepackage{mathrsfs}
\usepackage{amsbsy}
\usepackage{relsize}
\usepackage[colorlinks, linkcolor=red, citecolor=blue, urlcolor=blue, pagebackref, hypertexnames=false]{hyperref}
\usepackage[hyperpageref]{backref}

\usepackage{geometry}
\newcommand{\R}{\mathbb R}

\let\z=\zeta

\newtheorem{thm}{Theorem}[section]
\newtheorem{prop}[thm]{Proposition}
\newtheorem{lem}[thm]{Lemma}

\newtheorem{rem}[thm]{Remark}

\def\div{\mathop{\mathrm{div}}\nolimits}

\numberwithin{equation}{section}

\title[$p-$Laplacian heat equation]{Remarks on blow-up phenomena in $p-$Laplacian heat equation with inhomogeneous nonlinearity}

\author[E. A. Alzahrani and M. Majdoub]{ Eadah Ahmad Alzahrani and Mohamed Majdoub}

\address[E. A. Alzahrani]{Deapartment of Mathematics, College of Science, Imam Abdulrahman Bin Faisal University, P. O. Box 1982, Dammam, Saudi Arabia \& Basic and Applied Scientific Research Center, Imam Abdulrahman Bin Faisal University, P.O. Box 1982, 31441, Dammam, Saudi Arabia}
\email{\sl ealzahrani@iau.edu.sa}
\address[M. Majdoub]{Deapartment of Mathematics, College of Science, Imam Abdulrahman Bin Faisal University, P. O. Box 1982, Dammam, Saudi Arabia \& Basic and Applied Scientific Research Center, Imam Abdulrahman Bin Faisal University, P.O. Box 1982, 31441, Dammam, Saudi Arabia}
\email{\sl mmajdoub@iau.edu.sa}

\begin{document}
\begin{abstract}
We investigate the  $p-$Laplace heat equation $u_t-\Delta_p u=\z(t)f(u)$ on a bounded smooth domain $\Omega\subset\R^N$. Using differential inequalities arguments, we prove blow-up results under suitable conditions on $\z,\, f$, and the initial data $u_0$. We also give an upper bound for the blow-up time in each case.
\end{abstract}





\maketitle
\date{today}


\section{Introduction}
\setcounter{equation}{0}
In the past decade a strong interest in the phenomenon of blow-up of solutions to various classes of nonlinear parabolic problems has been assiduously investigated. We refer the reader to the books \cite{Straug, QS} as well as to the survey paper \cite{BB}. Problems with constant coefficients were investigated  in \cite{PS}, and problems with time-dependent coefficients under homogeneous Dirichlet boundary conditions were treated in \cite{PaPh2013}. See also \cite{PP} for a related system. The question of blow-up for nonnegative classical solutions of $p-$Laplacian heat equations with various boundary conditions has attracted considerable attention in the mathematical community in recent years. See for instance \cite{Ding2016, LW2008, Mess, LLY}.

There are two effective techniques which has been employed to prove non-existence of global solutions: the concavity method (\cite{Lev}) and the eigenfunction method (\cite{Kap}). The later one was firstly used for bounded domains but can be adapted to the whole space $\R^N$. The concavity method and its variants were used in the study of many nonlinear evolution partial differential equations (see e.g. \cite{GP1, GP2, PuSe}).

In the present paper, we investigate the blow-up phenomena of solutions to the following nonlinear $p-$Laplacian heat equation:
\begin{equation}
\label{main}
\left\{
\begin{matrix}
u_t-\Delta_p u=\z(t)f(u),\quad x\in \Omega,\quad t>0,\\
u(t,x)=0,\quad x\in \partial\Omega,\quad t>0,\\
u(0,x)= u_0(x),\quad x\in\Omega,\\
\end{matrix}
\right.
\end{equation}
where $\Delta_p u:= \div\left(|\nabla u|^{p-2}\nabla u\right)$ is the $p-$Laplacian operator, $p\geq 2$, $\Omega$ is a bounded sufficiently smooth domain in $\R^N$, $\z(t)$ is a nonnegative continuous function. The nonlinearity $f(u)$ is assumed to be continuous with $f(0)=0$. More specific assumptions on $f$, $\z$ and $u_0$ will be made later.

The case of $p=2$ has been studied in \cite{PS} for $\z(t)\equiv 1$, and in \cite{PaPh2013} for $\z$ being a non-constant function of $t$. Concerning the case $p>2$, Messaoudi \cite{Mess} proved the blow-up of solutions with vanishing initial energy when $\z(t)\equiv 1$. See also \cite{LW2008} and references therein. Recently, a $p-$Laplacian heat equations with nonlinear boundary conditions and time-dependent coefficient was investigated in \cite{Ding2016}. This note may be regarded as a complement, and in some sense an improvement, of \cite{PaPh2013, Mess}.\\

Let us now precise the assumptions on $f$ and $\z$. If $p=2$, we suppose either
\begin{equation}
\label{f1}
f\in C^1(\R) \quad \mbox{is convex with}\quad f(0)=0;
\end{equation}
\begin{equation}
\label{f2}
\exists\;\;\; \lambda>0\quad\mbox{such that}\quad f(s)>0\quad\mbox{for all}\quad s\geq \lambda;
\end{equation}
\begin{equation}
\label{f3}
\displaystyle\int\limits^\infty\,\frac{ds}{f(s)}<\infty;
\end{equation}
\begin{equation}
\label{z1}
\inf_{t\geq 0}\,\left(\int\limits_0^t\,(\z(s)-1)\,ds\right):=m\in(-\infty,0].
\end{equation}
or
\begin{equation}
\label{ff1}
s f(s)\geq (2+\epsilon)F(s)\geq C_0 |s|^\alpha,
\end{equation}
for some constants $\epsilon, C_0>0$, $\alpha>2$, and
\begin{equation}
\label{z2}
\z \in C^1([0,\infty))\quad \mbox{with}\quad \z(0)>0\;\;\;\mbox{and}\;\;\;\z'\geq 0.
\end{equation}
 Here $F(s)=\int\limits_0^s\,f(\tau)\,d\tau$.

Our first main result concerns the case $p=2$ and reads as follows.
\begin{thm}
\label{Th1}
Suppose that assumptions \eqref{f1}-\eqref{f2}-\eqref{f3}-\eqref{z1} are fulfilled. Let $0\leq u_0\in L^\infty(\Omega)$ satisfy $\int\limits_{\Omega}\,u_0\,\phi_1$ is large enough. Then the solution $u(t,x)$ of problem \eqref{main} blows up in finite time.
\end{thm}
\begin{rem}\quad\\
\vspace{-0.5cm}
\begin{itemize}
\item[(i)] The function $\phi_1$ stands for the eigenfunction of the Dirichlet-Laplace operator associated to the first eigenvalue $\lambda_1>0$, that is
$$
\Delta\phi_1=-\lambda_1\phi_1,\quad \phi_1>0,\;\;x\in\Omega; \quad \phi_1=0,\;\; x\in\partial\Omega, \;\;\int\limits_{\Omega}\phi_1=1.
$$
\item[(ii)] The assumptions \eqref{f1}-\eqref{f2}-\eqref{f3}-\eqref{z1} on $f$ and $\z$ cover the example
\begin{equation}
\label{Examp}
f(u)={\rm e}^{u}-1 \quad\mbox{and}\quad \z(t)={\rm e}^{t^2}.
\end{equation}
Note that this example is not studied in \cite{PaPh2013}, and Theorem \ref{Th1} can be seen as an improvement of Theorem 1 of \cite{PaPh2013}.
\item[(iii)] As it will be clear in the proof below, an upper bound of the maximal time of existence is given by
\begin{equation}
\label{Tmax1}
T^*=-m+2\int\limits_{y_0}^\infty\,\frac{ds}{f(s)},
\end{equation}
where $m$ is as in \eqref{z1} and $y_0={\rm e}^{m\lambda_1}\,\int\limits_{\Omega}\,u_0\phi_1$.
\item[(iv)] The conclusion of Theorem \ref{Th1} remains valid for $\Omega=\R^N$ if we replace $\phi_1$ by $\varphi(x)=\pi^{-N/2}\,{\rm e}^{-|x|^2}$.
\end{itemize}

\end{rem}

In order to state our next result (again for $p=2$), we introduce the energy functional
\begin{equation}
\label{energy}
E(u(t)):=\frac{1}{2}\int\limits_{\Omega}\,|\nabla u(t,x)|^2\,dx-\z(t)\int\limits_{\Omega}\,F(u(t,x))\,dx.
\end{equation}
Using \eqref{z2}, we see that $t\longmapsto E(u(t))$ is nonincreasing along any solution of \eqref{main}. This leads to the following.
\begin{thm}
\label{Th2}
Suppose that assumptions \eqref{ff1}-\eqref{z2} are fulfilled. Assume that either $E(u_0)\leq 0$ or $E(u_0)> 0$ and $\|u_0\|_{2}$ is large enough. Then the corresponding solution $u(t,x)$ blows up in finite time.
\end{thm}
\begin{rem}
An upper bound for the blow-up time is given by
\begin{eqnarray}
\label{Tmax2}
 T^*&=&\; \left\{
\begin{array}{cllll}\frac{(2+\epsilon)|\Omega|^{{\alpha/2}-1}\|u_0\|_2^{2-\alpha}}{\epsilon\z(0)C_0(\alpha-2)}
\quad&\mbox{if}&\quad
E(u_0)\leq 0,\\\\ \displaystyle\int\limits_{{\|u_0\|_2^2/2}}^\infty\,\frac{dz}{A z^{\alpha/2}-2E(u_0)} \quad
&\mbox{if}&\quad E(u_0)>0,
\end{array}
\right.
\end{eqnarray}
where
$$
A=\frac{2^{\alpha/2}C_0\epsilon \z(0)}{(2+\epsilon)|\Omega|^{{\alpha/2}-1}}.
$$
\end{rem}

We turn now to the case $p>2$. In \cite{Junning}, the author studied \eqref{main} when $\z(t)\equiv 1$. He established:
 \begin{itemize}
 \item  local existence when $f\in C^1(\R)$;
 \item global existence when $u f(u)\lesssim |u|^q$ for some $q<p$;
 \item nonglobal existence under the condition
 \begin{equation}
 \label{Nong}
 \frac{1}{p}\int\limits_{\Omega}\,|\nabla u_0)|^p\,dx-\int\limits_{\Omega}\,F(u_0)\,dx<0.
 \end{equation}
\end{itemize}
Later on Messaoudi \cite{Mess} improved the condition \eqref{Nong} by showing that blow-up can be obtained for vanishing initial energy. Note that by adapting the arguments used in \cite{Junning}, we can show a local existence result as stated below.
\begin{thm}
\label{Locex}
Suppose $\z\in C([0,\infty])$ and $f\in C(\R)$ satisfy $|f|\leq g$ for some $C^1-$function $g$. Then for any $u_0\in L^\infty(\Omega)\cap W_0^{1,p}(\Omega)$, the problem \eqref{main} has a local solution
$$
u\in L^\infty\left((0,T)\times\Omega\right)\cap L^p((0,T); W_0^{1,p}(\Omega)),\;\;\; u_t\in L^2\left((0,T)\times\Omega\right).
$$
\end{thm}

The energy of a solution $u$ is
\begin{equation}
\label{energy}
{\mathbf E}_p(u(t))={1\over p}\int\limits_{\Omega}\,|\nabla u(t,x)|^pdx-\z(t)\int\limits_{\Omega}\,F((u(t,x))dx\,.
\end{equation}
We also define the following set of initial data
\begin{equation}
\label{InDa}
{\mathcal E}=\Big\{\, u_0\in L^\infty(\Omega)\cap {W}^{1,p}_0(\Omega);\;\;\; u_0\not\equiv 0\;\;\;\mbox{and}\;\;\; {\mathbf E}_p(u_0)\leq 0\,\Big\}\,.
\end{equation}
Our main result concerning $p>2$ ca be stated as follows.
\begin{thm}
\label{Th3}
Suppose that assumption \eqref{z2} is fulfilled. Let $f\in C(\R)$ satisfy $|f|\leq g$ for some $C^1-$function $g$ and
\begin{equation}
\label{fF}
0\leq \kappa\,F(u)\leq uf(u),\quad \kappa>p>2.
\end{equation}
Then for any $u_0\in {\mathcal E}$ the solution $u(t,x)$ of \eqref{main} given in Theorem \ref{Locex} blows up in finite time.
\end{thm}
\begin{rem}
Although the proof uses the Poincar\'e inequality in a crucial way, we believe that a similar result can be obtained for $\Omega=\R^N$. This will be investigated in a forthcoming paper.
\end{rem}

We stress that the set ${\mathcal E}$ is non empty as it is shown in the following proposition.
\begin{prop}
\label{nonempty}
Suppose that assumption \eqref{fF} is fulfilled and $\z(0)>0$. Then ${\mathcal E}\neq\emptyset$.
\end{prop}

\section{Proofs}
This section is devoted to the proof of Theorems \ref{Th1}-\ref{Th2}-\ref{Th3} as well as Proposition \ref{nonempty}.
\subsection{Proof of Theorem \ref{Th1}}
The main idea in the proof is to define a suitable auxiliary function $y(t)$ and obtain a differential inequality leading to the blow-up. Define the function $y(t)$ as
\begin{equation}
\label{y}
y(t)=a(t)\int\limits_{\Omega}\,u(t,x)\,\phi_1(x)\,dx,
\end{equation}
where
\begin{equation}
\label{a}
a(t)={\rm e}^{\lambda_1(m-\Theta(t))},
\end{equation}
and

\begin{equation}
\label{Theta}
\Theta(t)=\int\limits_0^t\,\left(\z(s)-1\right)\,ds.
\end{equation}
We compute
\begin{eqnarray*}
y'(t)&=&\frac{a'(t)}{a(t)}\,y(t)-\lambda_1\,y(t)+a(t)\z(t)\int\limits_{\Omega}\,f(u(t,x))\phi_1(x)\,dx\\
&=&-\lambda_1\z(t)y(t)+a(t)\z(t)\int\limits_{\Omega}\,f(u(t,x))\phi_1(x)\,dx,
\end{eqnarray*}
where we have used $\frac{a'}{a}-\lambda_1=-\lambda_1\z$. By using \eqref{f1} and the fact that $0\leq a\leq 1$, we easily arrive at
\begin{equation}
\label{DI1}
y'(t)\geq\,\z(t)\left(-\lambda_1 y(t)+f(y(t))\right).
\end{equation}
Since $f$ is convex and due to \eqref{f3}, there exists a constant $C\geq \lambda$ such that $f(s)\geq 2\lambda_1 s$ for all $s\geq C$. Suppose $y(0)>C$. It follows from \eqref{DI1} that, as long as $u$ exists, $y(t)\geq C$. Therefore
$$
y(t)\geq \frac{\z(t)}{2}\,f(y(t)).
$$
Hence
$$
\frac{t+m}{2}\leq \frac{1}{2}\,\int\limits_0^t\,\z(s)\,ds\leq \int\limits_{y(0)}^\infty\,\frac{ds}{f(s)}<\infty.
$$
This means that the solution $u$ cannot exist globally and leads to the upper bound given by \eqref{Tmax1}.
\subsection{Proof of Theorem \ref{Th2}}
Let $y(t)$ be the auxiliary function defined as follows
$$
y(t)=\frac{1}{2}\int\limits_{\Omega}\,u^2(t,x)\,dx.
$$
We compute
\begin{eqnarray*}
y'(t)&=&\int\limits_{\Omega}\,u\left(\Delta u+\z(t)f(u)\right)\,dx\\
&=&-\int\limits_{\Omega}\,|\nabla u|^2\,dx+\z(t)\int\limits_{\Omega}\, u f(u)\,dx\\
&=&-2 E(u(t))+\z(t)\int\limits_{\Omega}\,\left(u f(u)-F(u)\right)\,dx,
\end{eqnarray*}
where $E(u(t))$ is given by \eqref{energy}. Taking advantage of \eqref{ff1}, we obtain that
\begin{equation}
\label{DI2}
y'(t)\geq -2 E(u(t))+\frac{\epsilon C_0}{2+\epsilon}\z(t)\,\int\limits_{\Omega}\,|u|^{\alpha}\,dx.
\end{equation}
Moreover, we compute
\begin{equation}
\label{En}
E'(u(t))=-\int\limits_{\Omega}\,u_t^2\,dx-\z'(t)\int\limits_{\Omega}\,F(u)\,dx\leq 0,
\end{equation}
thanks to \eqref{z2}. It then follows that $E(u(t))$ is non-decreasing in $t$ so that we have
\begin{equation}
\label{Ene}
E(u(t))\leq E(u(0))=E(u_0),\quad t\geq 0.
\end{equation}
From \eqref{DI2}, \eqref{Ene}, and the H\"older inequality, we find that
\begin{equation}
\label{DI3}
y'(t)\geq -2E(u_0)+\frac{\epsilon\z(0)C_0 2^{\alpha/2}}{(2+\epsilon) |\Omega|^{{\alpha/2}-1}}\,y(t)^{\alpha/2}.
\end{equation}
To conclude the proof we use the following result.
\begin{lem}
\label{Lem}
Let $y :[0,T)\to [0,\infty)$ be a $C^1-$function satisfying
\begin{equation}
\label{DI4}
y'(t)\geq -C_1+C_2\, y(t)^p,
\end{equation}
for some constants $C_1\in\R,\, C_2>0,\, p>1$. Then
\begin{eqnarray}
\label{Tmax2}
 T&\leq &\; \left\{
\begin{array}{cllll}\frac{y^{1-p}(0)}{C_2(p-1)}
\quad&\mbox{if}&\quad
C_1\leq 0,\\\\ \displaystyle\int\limits_{y(0)}^\infty\,\frac{dz}{C_2\, z^{p}-C_1} \quad
&\mbox{if}&\quad C_1>0.
\end{array}
\right.
\end{eqnarray}
\end{lem}
\subsection{Proof of Theorem \ref{Th3}}
We define
\begin{equation}
\label{H}
H(t)=\z(t)\int\limits_{\Omega}\,F((u(t,x))dx-{1\over p}\int\limits_{\Omega}\,|\nabla u(t,x)|^pdx,
\end{equation}
and
\begin{equation}
\label{L}
L(t)=\frac{1}{2}\|u(t)\|_2^2.
\end{equation}
By using \eqref{main}, we obtain that
\begin{eqnarray*}
H'(t) &=&\int\limits_{\Omega }u_{t}^{2}(t,x))\,dx+\z'(t)\,\int\limits_{\Omega
}F(u(t,x))\,dx \\
&=&\frac{\z'(t)}{\z(t)}\,H(t)+\int\limits_{\Omega }\,u_{t}^{2}(t,x))\,dx+
\frac{\z'(t)}{p\z(t)}\int\limits_{\Omega }|\nabla u(t,x)|^{p}\,dx \\
&\geq &\frac{\z'(t)}{\z(t)}\,H(t).
\end{eqnarray*}
Hence $H(t)\geq H(0)\geq 0$, by virtue of \eqref{z2}.

Recalling \eqref{main}, \eqref{H}, and \eqref{fF}, we compute
\begin{eqnarray*}
L'(t) &=&-\int\limits_{\Omega }|\nabla u(t,x)|^p\,dx+\z(t)\,
\int\limits_{\Omega }u(t,x)\,f(u(t,x))\,dx \\
&\geq &-\int\limits_{\Omega }|\nabla u(t,x)|^{p}dx+\kappa\,\z(t)\int\limits_{\Omega }\,F(u(t,x))dx\\
&\geq& \kappa\,H(t)+\left(\frac{\kappa}{p}-1\right)\,\int\limits_{\Omega}\,|\nabla u(t,x)|^p\,dx\\
&\geq& \left(\frac{\kappa}{p}-1\right)\,\int\limits_{\Omega}\,|\nabla u(t,x)|^p\,dx.
\end{eqnarray*}
Applying H\"older inequality and then Poincar\'e inequality yields
$$
L(t)\leq |\Omega|^{1-{2/p}}\, \left(\int\limits_{\Omega}\, |u(t,x)|^p\,dx\right)^{2/p}\leq C\,\left(\int\limits_{\Omega}\, |\nabla u(t,x)|^p\,dx\right)^{2/p},
$$
where $C>0$ is a constant depending only on $\Omega$ and $p$. Hence
\begin{equation}
\label{DI5}
L'(t)\geq\, \frac{\kappa-p}{p\,C^{p/2}}\, L^{p/2}(t).
\end{equation}
Integrating the  differential inequality \eqref{DI5} leads to
$$
t\leq \frac{2p C^{p/2} L^{1-p/2}(0)}{(p-2)(\kappa-p)}<\infty.
$$
Therefore $u$ blows up at a finite time $T^*\leq \frac{2p C^{p/2} L^{1-p/2}(0)}{(p-2)(\kappa-p)}$.
\subsection{Proof of Proposition\ref{nonempty}}
Recalling \eqref{fF}, we obtain that
\begin{equation}
\label{growth}
F(u)\geq C\, u^{\kappa}, \;\;\;\mbox{for all}\;\;\; u\geq 1,
\end{equation}
for some  constant $C>0$. Let $K\subset \Omega$ be a compact nonempty subset of $\Omega$. Fix a smooth cut-of  function $\phi\in C^\infty(\Omega)$ such that
$$
\phi(x)=1 \quad\mbox{for}\quad x\in K.
$$
 We look for  initial data $u_0=\lambda \phi$ where $\lambda>0$ to be chosen later. Clearly $u_0 \in L^\infty(\Omega)\cap {W}^{1,p}_0(\Omega)$, and for $\lambda\geq 1$ we have using \eqref{growth}
\begin{eqnarray*}
{\mathbf E}_p(u_0)&=&\frac{1}{p}\int\limits_{\Omega}\, |\nabla u_0|^p-\z(0)\int\limits_{\Omega}\, F(u_0),\\
&=&\frac{\lambda^p}{p}\int\limits_{\Omega}\, |\nabla\phi|^p-\z(0)\int\limits_{K}\, F(\lambda)-\z(0)\int\limits_{\Omega\backslash K}\, F(u_0),\\
&\leq&\frac{\lambda^p}{p}\int\limits_{\Omega}\, |\nabla\phi|^p - \tilde{C} \lambda^{\kappa},
\end{eqnarray*}
for some constant $\tilde{C}>0$. Since $\frac{\lambda^p}{p}\int\limits_{\Omega}\, |\nabla\phi|^p - \tilde{C} \lambda^{\kappa}\leq 0$ for $\lambda\geq \Big(\frac{\|\nabla\phi\|_p^p}{p\tilde{C}}\Big)^{1/(\kappa-p)}$, we deduce that $u_0\in \mathcal{E}$ for $\lambda$ large enough. This finishes the proof of Proposition \ref{nonempty}.


\end{document}